\documentclass[12pt,a4paper,twoside,reqno, final]{amsart}

\usepackage[notref,notcite]{showkeys}

\usepackage[all,cmtip]{xy}
\usepackage{amsmath,amscd}
\usepackage{amssymb,amsfonts}
\usepackage[driver=pdftex,margin=3cm,heightrounded=true,centering]{geometry}
\usepackage{url}
\usepackage[colorlinks=true,linkcolor=blue,citecolor=blue]{hyperref}
\usepackage{slashed}
\usepackage[T1]{fontenc}

\usepackage{graphicx}

\usepackage{a4wide}

\usepackage{thmtools}
\declaretheorem[style=theorem,name={Theorem}]{theoremletter}

\newtheorem{introcorollary}[theoremletter]{Corollary}

\usepackage{amsthm}

\theoremstyle{plain}

\newtheorem{theorem}{Theorem}[section]

\newtheorem{prop}[theorem]{Proposition}
\newtheorem{lemma}[theorem]{Lemma}
\newtheorem{cor}[theorem]{Corollary}

\theoremstyle{definition}
\newtheorem{dfn}[theorem]{Definition}

\theoremstyle{remark} 

\newtheorem{remark}[theorem]{Remark}

\theoremstyle{plain}

\usepackage{amscd}

\numberwithin{equation}{section}

\newcommand{\alpheqn}[1][\relax]{
     \refstepcounter{equation}
     \if#1\relax \relax
       \else \label{#1}
     \fi  
     \setcounter{saveeqn}{\value{equation}}%
    \setcounter{equation}{0}%
    \renewcommand{\theequation}{\thealphequation}}
\newcommand{\reseteqn}{\setcounter{equation}{\value{saveeqn}}%
     \renewcommand{\theequation}{\thearabicequation}}

\providecommand{\mathscr}{\mathcal} 
\IfFileExists{mathrsfs.sty}{
\usepackage{mathrsfs}}         
   {
     \IfFileExists{eucal.sty}{
        \usepackage[mathscr]{eucal}} 
   {
   }
}

\newcommand{\Lip}{\operatorname{Lip}}

\newcommand{\sa}{{\operatorname{sa}}}

\newcommand{\vertiii}[1]{{\left\vert\kern-0.25ex\left\vert\kern-0.25ex\left\vert #1 
    \right\vert\kern-0.25ex\right\vert\kern-0.25ex\right\vert}}
\newcommand{\Bvert}[1]{{\Big\vert\kern-0.25ex\Big\vert\kern-0.25ex\Big\vert #1 
    \Big\vert\kern-0.25ex\Big\vert\kern-0.25ex\Big\vert}}
\newcommand{\bvert}[1]{{\big\vert\kern-0.25ex\big\vert\kern-0.25ex\big\vert #1 
    \big\vert\kern-0.25ex\big\vert\kern-0.25ex\big\vert}}
\newcommand{\nvert}[1]{{\vert\kern-0.25ex\vert\kern-0.25ex\vert #1 
    \vert\kern-0.25ex\vert\kern-0.25ex\vert}}

\renewcommand{\leq}{\leqslant}
\renewcommand{\geq}{\geqslant}

\newcommand{\cd}{\cdot}

\newcommand{\ot}{\otimes}
\newcommand{\hot}{\widehat \otimes}

\newcommand{\op}{\oplus}

\newcommand{\ci}{\circ}

\newcommand{\ti}{\times}
\newcommand{\nn}{\mathbb{N}}
\newcommand{\zz}{\mathbb{Z}}

\newcommand{\rr}{\mathbb{R}}
\newcommand{\cc}{\mathbb{C}}

\newcommand{\be}{\beta}
\newcommand{\ga}{\gamma}

\newcommand{\de}{\delta}
\newcommand{\De}{\Delta}
\newcommand{\ep}{\varepsilon}

\newcommand{\La}{\Lambda}

\newcommand{\ze}{\zeta}
\newcommand{\pa}{\partial}

\newcommand{\ov}{\overline}
\newcommand{\C}[1]{\mathcal{#1}}
\newcommand{\G}[1]{\mathfrak{#1}}
\newcommand{\T}[1]{\textup{#1}}

\newcommand{\B}[1]{\mathbb{#1}}

\newcommand{\s}[1]{\mathscr{#1}}

\newcommand{\fork}[2]{\left\{ \begin{array}{#1} #2 \end{array} \right.} 

\newcommand{\ma}[2]{\left(\begin{array}{#1} #2 \end{array} \right)}

\newcommand{\su}{\subseteq}

\newcommand{\q}{\qquad}

\newcommand{\wit}{\widetilde}

\newcommand{\inn}[1]{\langle #1 \rangle}
\newcommand{\binn}[1]{\big\langle #1 \big\rangle}

\author{Konrad Aguilar}

\address{ Department of Mathematics, Pomona College, 610 N. College Ave., Claremont, CA 91711}
\email{konrad.aguilar@pomona.edu}

\author{Jens Kaad}
\address{Department of Mathematics and Computer Science,
The University of Southern Denmark,
Campusvej 55, DK-5230 Odense M,
Denmark}
\email{kaad@imada.sdu.dk}

\author{David Kyed}
\address{Department of Mathematics and Computer Science,
The University of Southern Denmark,
Campusvej 55, DK-5230 Odense M,
Denmark}
\email{dkyed@imada.sdu.dk}

\title{Polynomial approximation of quantum Lipschitz functions}



\begin{document}

\begin{abstract}
We prove an approximation result for Lipschitz functions on the quantum sphere $S_q^2$, from which we deduce that the two natural quantum metric structures on $S_q^2$ have quantum Gromov-Hausdorff distance zero. 
\end{abstract}

\maketitle

\tableofcontents


\section{Introduction}
The theory of operator algebras forms the core of many interesting non-commutative generalisations  of classical mathematical theories, 
including non-commutative topology \cite{Bla:K-theory},  free probability \cite{MiSp:Free}, quantum groups \cite{KlSc:QGR} and non-commutative geometry \cite{Con:NCG}.
Within this general paradigm,  Rieffel's theory of compact quantum metric spaces \cite{Rie:MSS, Rie:GHD} provides an elegant non-commutative counterpart to classical compact metric spaces. 
 The essential data defining such a compact quantum metric structure is given by a densely defined seminorm $L$ on a unital $C^*$-algebra $A$. The main requirement for $L$ is that the extended metric
\[
\rho_L(\mu, \nu):=\sup\{|\mu(a)-\nu(a)| \, \mid \, a\in \T{Dom}(L), L(a)\leq 1\} \in [0,\infty], \qquad \mu,\nu\in \C S(A),
\]
defines a genuine (i.e.~everywhere finite) metric on the state space $\C S(A)$ and that this metrises the weak$^*$-topology, in which case $L$ is referred to as a \emph{Lip-norm}.
The definition is influenced by Connes' non-commutative geometry \cite{Con:NCG}, since one naturally  obtains a seminorm  from  a unital spectral triple $(\C A, H, D)$ by setting
\begin{align}\label{eq:commutator-seminorm}
L_D(a):=\big\|\ov{[D, a]}\big\|, \quad a\in \C A.
\end{align}
For a given $C^*$-algebra $A$ and Lip-norm $L$, it may be possible to enlarge or reduce the domain of $L$ and thereby obtain different quantum metric structures on $A$. In the case of a unital spectral triple $(\C A,H,D)$, there is a natural maximal seminorm $L_D^{\max}$ defined by the formula \eqref{eq:commutator-seminorm} but with domain
\[
A^{\Lip}:=\big\{ a\in A \mid a(\T{Dom} (D))\subseteq \T{Dom}(D) \T{ and } [D,a] \T{ extends boundedly to $H$}   \big\}.
\]
Here $A$ denotes the norm closure of $\C A\subseteq \mathbb{B}(H)$. Experience with concrete examples shows that it may not be possible to recover $L_D^{\max}$ from $L_D$. In fact, the assignment $a \mapsto \ov{[D,a]}$ yields a closable derivation $\pa \colon \C A \to \B B(H)$ and the closure of $\pa$ yields an intermediate algebra $\C A \su A^1 \su A^{\Lip}$ which is, in general, different from both $\C A$ and $A^{\Lip}$. Whereas all the relevant information regarding the extension $L_D^1 \colon A^1 \to [0,\infty)$ can be obtained from $L_D \colon \C A \to [0,\infty)$ by approximation arguments, this is not the case for $L_D^{\max}$. In fact, the analysis of $L_D^{\max}$ requires different methods relying more on von Neumann algebraic techniques than $C^*$-algebraic techniques. In the classical case of a compact spin manifold, the three different domains (coming from the Dirac operator) would be smooth functions (or an appropriate algebra of polynomials), $C^1$-functions and Lipschitz functions. Notice in this respect that the derivative of a Lipschitz function makes sense but only in the von Neumann algebraic context of (equivalence classes of) bounded measurable functions.

One of the main virtues of the theory of compact quantum metric spaces is that it allows for an analogue of the classical Gromov-Hausdorff distance \cite{edwards-GH-paper, gromov-groups-of-polynomial-growth-and-expanding-maps}, known as the \emph{quantum Gromov-Hausdorff distance} \cite{Rie:GHD} and denoted $\T{dist}_{\T{Q}}$, see Section \ref{subsec:cqms} for the definition.  This allows one to study the class of compact quantum metric spaces from an analytical point of view, and ask questions pertaining to continuity and convergence of families of compact quantum metric spaces, see e.g.~\cite{aguilar:thesis,  kaad-kyed,  Lat:AQQ, LatPack:Solenoids,  Rie:GHD,  Rie:MSG} for examples of this.

The main focus in the present paper is the Podle{\'s} sphere $S_q^2$ \cite{Pod:QS}, which forms the base of a spectral triple for the D\k{a}browski-Sitarz Dirac operator $D_q$ \cite{DaSi:DSP}, whose associated seminorm $L_{D_q}^{\max}$ turns $C(S_q^2)$ into a compact quantum metric space, as proven in \cite{AgKa:PSM}. However, the natural point of departure when studying $S_q^2$ is actually the associated coordinate algebra $\C O (S_q^2)$ which is a subalgebra of the \emph{Lipschitz algebra} $C^{\Lip}(S_q^2):=\T{Dom}(L_{D_q}^{\max})$, and one may therefore restrict $L_{D_q}^{\max}$ to $\C O(S_q^2)$ to obtain another Lip-norm $L_{D_q}$. In \cite{AKK:PodCon} we undertook a detailed study of this Lip-norm, and proved that the family of compact quantum metric spaces $(C(S_q^2), L_{D_q})_{q\in (0,1]}$ varies continuously in $q$ with respect to the quantum Gromov-Hausdorff distance --  thus in particular showing that the quantised 2-spheres $S_q^2$ converge to the classical round 2-sphere $S^2$ as $q$ tends to 1. Although the coordinate algebra $\C O(S_q^2)$ is a very natural domain for the Lip-norm when approaching the theory of $q$-deformed spaces from a Hopf-algebraic  angle, the Lipschitz algebra $C^{\Lip}(S_q^2)$ is the more natural domain from the point of view of non-commutative geometry. So the question remaining is if the convergence results from \cite{AKK:PodCon} hold true also when $L_{D_q}$ is replaced with $L_{D_q}^{\max}$. The main point of the present paper is to answer this in the affirmative, by proving the following: 

\begin{theoremletter}\label{mainthm:A}
It holds that $\T{dist}_{\operatorname{Q}} \big( (C(S_q^2), L_{D_q}); (C(S_q^2), L_{D_q}^{\max})  \big)=0$.
\end{theoremletter}
As explained earlier in this introduction we also have an intermediate algebra $C^1(S_q^2)$ obtained by taking the closure of the derivation coming from $D_q$ and it is straightforward to show that $\T{dist}_{\operatorname{Q}} \big( (C(S_q^2), L_{D_q}); (C(S_q^2), L_{D_q}^1)  \big) = 0$ (for more general considerations along these lines, see \cite[Proposition 7.1]{Rie:GHD}). The fact that this distance is equal to 0 is however of little help when approaching the more general result in Theorem \ref{mainthm:A}. The algebra $C^1(S_q^2)$ is indeed substantially different from $C^{\Lip}(S_q^2)$ and the relevant information regarding the seminorm $L_{D_q}^{\max}$ cannot be retrieved from the seminorm $L_{D_q}$ by standard approximation techniques. The deep gap between the two seminorms in question is perhaps best illustrated by recalling the difference between the $C^*$-algebra $C(S_q^2)$ and the enveloping von Neumann algebra $L^\infty(S_q^2)$: Similarly to the classical situation, derivatives of elements in the Lipschitz algebra $C^{\Lip}(S_q^2)$ are, in general, not elements in $C(S_q^2)$ but can only be described using the von Neumann algebraic framework of $L^\infty(S_q^2)$, see Lemma \ref{lem:delta-in-L-infty} for a more detailed statement.\\

The result in Theorem \ref{mainthm:A} is a consequence of a careful analysis of the quantum Berezin transform introduced in \cite{AKK:PodCon}. In fact, our analysis yields the following approximation result, which is of interest in its own right:

\begin{theoremletter}\label{mainthm:B}
For every $q\in (0,1]$, every $\ep>0$ and every $x\in C^{\Lip}(S_q^2)$ there exists $y\in \C O(S_q^2)$ with $L_{D_q}(y)\leq L_{D_q}^{\max}(x)$ and $\|x-y\|\leq \ep \cd L_{D_q}^{\max}(x)$.
\end{theoremletter}

For $q=1$, one has $C(S_1^2)=C(S^2)$ and $\C O(S^2)$ is generated by the three coordinate functions $x_1,x_2,x_3$ arising from the natural embedding of $S^2$ into $\rr^3$. Moreover,  $L_{D_1}^{\max}$ recovers the Lipschitz constant $L_{d_{S^2}}(f)$ of a function $f$ on $S^2$ with respect to the usual round metric $d_{S^2}$,  and Theorem \ref{mainthm:B}  therefore in particular includes an approximation result for classical functions\footnote{We suspect that this approximation result for classical functions is well-known but we were unable to find it in the litterature.}:

\begin{introcorollary}
For every Lipschitz function $f$ on $S^2$ and every $\ep>0$ there exists $p\in \C O(S^2)$ with $L_{d_{S^2}}(p)\leq L_{d_{S^2}}(f)$ and $\| f-p \|_\infty< \ep \cdot  L_{d_{S^2}}(f)$.
\end{introcorollary}

We finally record that the combination of the main convergence result from \cite{AKK:PodCon} and the present Theorem \ref{mainthm:A} yield the following:

\begin{introcorollary}
For every $q_0 \in (0,1]$ we have that 
\[
\displaystyle\lim_{q \to q_0} \T{dist}_{\operatorname{Q}} \big( (C(S_q^2), L_{D_q}^{\max}); (C(S_{q_0}^2), L_{D_{q_0}}^{\max})  \big)=0 .
\]
\end{introcorollary}

\paragraph{\textbf{Acknowledgements.}}
The authors gratefully acknowledge the financial support from  the Independent Research Fund Denmark through grant no.~9040-00107B  and 7014-00145B. They are furthermore grateful to the anonymous referee for a careful reading of the manuscript and for helpful suggestions and comments.\\

\paragraph{\textbf{Notation.}}
Throughout the text, all inner products are assumed linear in the second variable. The symbols $\ot$, $\ot_{\min}, \hot$ and $\bar{\ot}$ will denote algebraic tensor products, minimal tensor products of $C^*$-algebras, tensor products of Hilbert spaces and tensor products of von Neumann algebras, respectively.

\section{Preliminaries}
In this section we give the necessary preliminaries concerning compact quantum metric spaces, quantum $SU(2)$ and the Podle{\'s} sphere. We will align our notation with the one used in \cite{AKK:PodCon} where a much more in-depth introduction is given, and for this reason the present presentation will be kept relatively laconic.

\subsection{Compact quantum metric spaces}\label{subsec:cqms}
The theory of compact quantum metric spaces was initiated by Rieffel around year 2000 \cite{Rie:MSS, Rie:GHD}, and the first approach took the theory of order unit spaces as its point of departure. Since then, several variations of the theory have emerged \cite{Li:CQG, Li:GH-dist, Rie:MSS}, and we here take an operator system approach aligning with the recent developments in  \cite{walter-connes:truncations, walter:GH-convergence}. Recall that a \emph{concrete operator system} is a closed subspace $X$ of a unital $C^*$-algebra such that $X$ is stable under the involution and contains the unit. For a $*$-invariant unital subspace $\C X \su X$ we let $\C X_{\sa} := \{ x \in \C X \mid x = x^* \}$ denote the real vector space of selfadjoint elements in $\C X$. Every operator system $X$ has a state space $\C S(X)$ consisting of all the positive linear functionals mapping the unit to 1, and this leads to the following definition of a compact quantum metric space:

\begin{dfn}\label{def:cqms}
A \emph{compact quantum metric space} is a concrete operator system $X$ equipped with a densely defined seminorm $L\colon \T{Dom}(L) \to [0,\infty)$ satisfying that:
\begin{itemize}
\item[(i)] $\cc \cd 1 \su \T{Dom}(L)$ and $L(x)=0$ if and only if $x\in \mathbb{C}\cdot 1$;
\item[(ii)] $x^* \in \T{Dom}(L)$ and $L(x)=L(x^*)$ for all $x\in \T{Dom}(L)$;
\item[(iii)] The function $\rho_L(\mu,\nu):=\sup\{|\mu(x)-\nu(x)| \mid x \in \T{Dom}(L) \T{ and } L(x)\leq 1\}$ equips the state space $\C S(X)$ with a metric which metrises the weak$^*$-topology.
\end{itemize}
A densely defined seminorm $L$ satisfying (i)-(iii) is called a \emph{Lip-norm} and the corresponding metric $\rho_L$ is referred to as the \emph{Monge-Kantorovi\v{c} metric}.
\end{dfn}

\begin{remark}
For a compact quantum metric space $(X,L)$ one often extends $L$ to all of $X$ by setting it equal to infinity outside of $\T{Dom}(L)$, which of course captures the same information. This is for instance the approach taken in \cite{AKK:PodCon}, but since the domain of definition is particularly important in the present paper we keep the seminorms on their finite domains to avoid confusion.
\end{remark}

As already mentioned in the introduction, unital spectral triples provide a natural source of examples of compact quantum metric spaces, but to better understand the origin of Rieffel's definitions, it is illuminating to also briefly discuss the commutative case. Starting out with a compact metric space $(X,d)$ and forming the associated $C^*$-algebra $C(X)$, one can consider the subalgebra $C^{\Lip}(X)$ consisting of Lipschitz continuous functions. On $C^{\Lip}(X)$ the metric gives rise to a natural seminorm defined by
\[
L_d(f):=\sup\Big\{ \frac{|f(x)-f(y)|}{d(x,y)} \, \Big\vert \, x,y\in X, x\neq y\Big\},
\]
and it is well known that $L_d$ is a Lip-norm, and that the restriction of $\rho_{L_d}$ to $X\subseteq \C S(C(X))$ agrees with the original metric $d$ \cite{KaRu:FSE,KaRu:OSC}.\\

If $(X,L_X)$ and $(Y,L_Y)$  are two compact quantum metric spaces, then a Lip-norm $L$ on $X \oplus Y$ with domain $\T{Dom}(L_X)\oplus \T{Dom}(L_Y)$ is called \emph{admissible} if the two quotient seminorms, the restriction
\[
L \colon \T{Dom}(L_X)_{\sa} \op \T{Dom}(L_Y)_{\sa} \to [0,\infty)
\]
induces via the natural projections, agree with the restrictions $L_X \colon \T{Dom}(L_X)_\sa \to [0,\infty)$ and $L_Y \colon \T{Dom}(L_Y)_\sa \to [0,\infty)$. For such an admissible $L$, one obtains isometric embeddings
\[
(\C S(X), \rho_{L_X}) \hookrightarrow (\C S (X\oplus Y), \rho_L) \quad \text{ and }  \quad (\C S(Y), \rho_{L_Y}) \hookrightarrow (\C S (X\oplus Y), \rho_L) 
\]
and hence the \emph{Hausdorff distance} $\T{dist}_{\text{H}}^{\rho_L}\big(\C S(X), \C S(Y)\big)$ makes sense \cite{Hausdorff-grundzuge}. The \emph{quantum Gromov-Hausdorff distance} between $(X,L_X)$ and $(Y,L_Y)$ is then defined as
\[
\T{dist}_{\text{Q}} ((X,L_X); (Y,L_Y)):= \inf \left\{ \T{dist}_{\text{H}}^{\rho_L}(\C S(X), \C S(Y)) \mid L \text{  admissible}   \right\}.
\]
We underscore that this is simply a reformulation of Rieffel's original definition based on order unit spaces \cite{Rie:GHD}.
More precisely, putting $A:=  \{x\in \T{Dom}(L_X) \mid x=x^*\}$ and $B:=\{y\in \T{Dom}(L_Y) \mid y=y^*\}$ we obtain order unit compact quantum metric spaces and $\T{dist}_{\text{Q}}\big((X,L_X); (Y,L_Y)\big)=\T{dist}_{\text{Q}}\big((A,L_X\vert_A); (B,L_Y\vert_B)\big)$. In particular, it is important to note that $\T{dist}_{\text{Q}}\big((X,L_X); (Y,L_Y)\big) = 0$ is equivalent to the existence of an isometric, affine bijection from $(\C S(Y), \rho_{L_Y})$ to
$(\C S(X), \rho_{L_X})$, or, equivalently, a Lip-norm isometric order unit isomorphism at the level of (the closures of) $A$ and $B$, see \cite[Corollary 6.4 \& Theorem 7.9]{Rie:GHD}. Several more sophisticated notions of quantum distances have been proposed over the past 20 years  \cite{Ker:MQG, Lat:QGH, Li:CQG}, but in the present paper we will only be concerned with Rieffel's original version.


\subsection{Quantum $SU(2)$ and the Podle{\'s} sphere}
In this section we briefly introduce the main objects of study and fix the notation.  As in the previous section,  we align  our notation with that of  \cite{AKK:PodCon},  which also contains a more in-depth introduction to the material covered below. The general literature on  quantum groups  is vast, and we refer the reader to the monographs  \cite{KlSc:QGR} and \cite{Tim:Invitation} and references therein for the background theory.  A central role in the present paper is played by Woronowicz' quantum $SU_q(2)$ introduced in \cite{Wor:UAC}  which is defined via a universal unital $C^*$-algebra $C(SU_q(2))$ with generators $a$ and $b$ subject to the relations
\[
\begin{split}
& ba = q ab \quad \quad b^* a = q ab^* \quad \quad bb^* = b^* b \\
& a^* a + q^2 bb^* = 1 = aa^* + bb^* .
\end{split}
\]
This becomes a  $C^*$-algebraic compact quantum group \cite{wor-cp-qgrps} and we denote its comultiplication by $\Delta\colon C(SU_q(2))\to C(SU_q(2))\ot_{\min} C(SU_q(2))$ and the Haar state (see \cite[Theorem 1.3]{wor-cp-qgrps}) by $h$. Furthermore, we fix a complete set of irreducible unitary corepresentations $(u^n)_{n\in \nn_0}$, satisfying the additional technical conditions specified  in \cite[Section 2.1]{AKK:PodCon}.  In particular, the corepresentation $u^1$ is the so-called \emph{fundamental unitary corepresentation} given by the matrix
\[
u^1 = \begin{pmatrix} a^* & -qb \\ b^* & a \end{pmatrix},
\]
which is also denoted $u$ in the sequel. We denote the GNS-space associated with $h$ by $L^2(SU_q(2))$, the natural inclusion $C(SU_q(2))\subseteq L^2(SU_q(2))$ by $\Lambda$ and the associated GNS-representation $C(SU_q(2))\to \mathbb{B}(L^2(SU_q(2)))$ by $\rho$. 
The matrix coefficients $u_{ij}^n$ are linearly independent and form an orthogonal basis of $L^2(SU_q(2))$ (once included in this Hilbert space via $\La$). We apply the notation
\[
\zeta_{ij}^n:= h\big((u_{ij}^n)^*u_{ij}^n\big)^{-1/2} \cdot u_{ij}^n \in C(SU_q(2))
\]
for the corresponding normalised elements, which then have the property that $\{\La(\zeta_{ij}^n)\mid n\in \nn_0, 0\leq i,j\leq n\}$ is an orthonormal basis for $L^2(SU_q(2))$.

In addition to the $C^*$-algebraic picture, one may consider $SU_q(2)$ from an algebraic point of view, by restricting attention to the unital $*$-subalgebra $\C O(SU_q(2))$ generated by $a$ and $b$. This unital $*$-subalgebra can be given the structure of a unital Hopf $*$-algebra where the comultiplication is obtained by restriction of the comultiplication on $C(SU_q(2))$. We shall denote the counit and antipode of $\C O(SU_q(2))$ by $\epsilon$ and $S$ respectively, and note that $SU_q(2)$ is known to be co-amenable \cite[Theorem 2.12]{bedos-murphy-tuset}, meaning that the counit extends boundedly to a character on $C(SU_q(2))$. We note that
\begin{align*}
\C O (SU_q(2))=\T{span}_{\cc}\{ u_{ij}^n \mid n\in \nn_0, 0\leq i,j, \leq n\}.
\end{align*}
Also at the Lie algebra level it is possible to pass to the $q$-deformed level, and upon doing so one arrives at the \emph{quantum enveloping algebra} $ \C U_q( \G{su}(2))$ with generators $e,f$ and $k$ (see \cite[Chapter 3]{KlSc:QGR} for the precise definition). The quantum enveloping algebra is also a unital Hopf $*$-algebra which admits a non-degenerate dual pairing with $\C O(SU_q(2))$ denoted by $\inn{\cd, \cd}\colon  \C U_q( \G{su}(2)) \ti \C O (SU_q(2))\to \cc$. The dual pairing, in turn, gives rise to right and left actions of $ \C U_q( \G{su}(2))$ on $\C O (SU_q(2))$ (by linear endomorphisms) defined, respectively, by
\[
\de_\eta(x) := (\inn{\eta,\cd} \ot 1)\De(x) \qquad \T{ and } \qquad \pa_\eta(x) := (1 \ot \inn{\eta,\cd})\De(x), 
\]
for $x\in \C O(SU_q(2))$ and $\eta\in \C U_q( \G{su}(2))$. Among these operators, the following three play a key role in the next section:
\begin{align}\label{eq:delta-definitions}
\de_1:= q^{1/2} \de_e, \qquad \de_2:= q^{-1/2}\de_f \qquad \T{and} \qquad \de_3:=\frac{\de_k-\de_{k^{-1}}}{q-q^{-1}} .
\end{align}
The operator $\de_3$ is of course only well-defined by the above formula for $q\neq 1$, and for $q=1$ we put  $\de_3:=\frac12 \de_h$ where $h = [f,e]$ in the classical enveloping Lie algebra $\C U( \G{su}(2))$, see \cite{AKK:PodCon} for a more detailed discussion of this. The operators $\de_1,\de_2$ and $\de_3$ are all $\de_k$-twisted derivations (with $k = 1$ for $q = 1$), meaning that 
\begin{align}\label{eq:twisted-leibniz}
\de_i(xy)=\de_i(x)\de_k(y) +\de_{k^{-1}}(x)\de_i(y), \qquad x,y\in \C{O}(SU_q(2)), \quad i\in \{1,2,3\}.
\end{align}
Moreover, they are compatible with the adjoint operation in the sense that
\begin{align}\label{eq:starder}
\de_1(x^*) = - \de_2(x)^* \q \T{and} \q \de_3(x^*) = - \de_3(x)^*, \qquad x \in \C O(SU_q(2)) .
\end{align} 
For each $\eta \in \C U_q(\G{su}(2))$ it holds that $h\circ \de_\eta=h\circ \pa_\eta=\eta(1)\cdot h$ which follows directly from the bi-invariance of the Haar state; in particular 
\begin{align}\label{eq:h-circ-delta=0}
h\circ \de_1=h\circ \de_2=h\circ \de_3=0 \quad \T{ and } \quad h\circ \de_k=h,
\end{align} 
by the properties of the pairing and the comultiplication in $\C U_q( \G{su}(2))$. Both $\de_k$ and $\pa_k$ are algebra automorphisms of $\C O(SU_q(2))$ and the composition $\nu:=\de_{k^{-2}}\circ \pa_{k^{-2}}$ is the modular automorphism (see \cite[Chapter 4, Proposition 15]{KlSc:QGR}) for the Haar state, meaning that 
\begin{equation}\label{eq:twisttrace}
h(xy)=h(\nu(y) x ), \qquad y \in \C O(SU_q(2)) \T{ and } x \in C(SU_q(2)) .
\end{equation}


The classical Hopf fibration $SU(2) \to S^2$ shows that $C(S^2)$ may be viewed as the fixed point algebra $C(SU(2))^{S^1}$ for the induced circle action $S^1\curvearrowright C(SU(2))$. The $q$-deformed analogue also admits a natural circle action $S^1\overset{\sigma}{\curvearrowright} C(SU_q(2))$  (given on generators by $\sigma_z(a)=za$ and $\sigma_z(b)=zb$) and  the Podle{\'s} sphere $S_q^2$ is defined, implicitly, via the $C^*$-algebra $C(S_q^2):= C(SU_q(2))^{S^1}$. Concretely, $C(S_q^2)$ is generated by the elements $A:=b^*b$ and $B:=ab^*$, and the dense unital $*$-algebra generated by $A$ and $B$ is denoted $\C O(S_q^2) \su C(S_q^2)$ and referred to as the \emph{coordinate algebra}. This unital $*$-algebra can also be described in terms of the matrix units, as one has the identity
 \begin{align*}
\C O (S_q^2)=\T{span}_{\cc}\{ u_{in}^{2n} \mid n\in \nn_0, 0\leq i \leq 2n\}.
\end{align*}
The aforementioned circle action gives rise to spectral subspaces
\begin{align}\label{eq:def-of-spectral-subspaces}
\C A_n :=\{ x\in \C O(SU_q(2)) \mid \forall z\in S^1: \sigma_z(x)=z^n x\}, \qquad n\in \zz
\end{align}
which may also be described (see e.g.~\cite[Theorem 1.2]{Wor:UAC}) as
\begin{align*}\label{eq:spectral-subspaces-in-terms-of-pa}
\C A_n :=\{ x\in \C O(SU_q(2)) \mid \pa_k(x)=q^{n/2}x\}, \qquad q \neq 1.
\end{align*}
Note that it follows from this description that the modular automorphism $\nu$ preserves the spectral subspaces, a fact we will be using in the sequel without further mentioning. We denote by $H_-$ and $H_+$ the closure of $\Lambda(\C A_{-1})$ and $\Lambda(\C A_{1})$ in $L^2(SU_q(2))$, respectively, while the closure of $\Lambda(\C O(S_q^2))$ will be denoted $L^2(S_q^2)$. Clearly, each $\C A_n$ is a left module over $\C A_0=\C O(S_q^2)$, so upon restricting the GNS-representation, both $H_-$ and $H_+$ acquire an action of $C(S_q^2)$ and hence so does their direct sum $H_+\oplus H_-$. We denote the corresponding representation by
\[
\pi\colon C(S_q^2)\to \mathbb{B}(H_+\oplus H_-) .
\]
In \cite{DaSi:DSP}, D{\c a}browski and Sitarz provided a Dirac operator $D_q$ on $H_+\oplus H_-$ and they proved (among other things) that the coordinate algebra $\C O(S_q^2)$ then fits in an even unital spectral triple $(\C O(S_q^2), H_+ \op H_-, D_q)$. Concretely, $D_q$ is defined as follows: The endomorphisms $\pa_e$ and $\pa_f$ of $\C O(SU_q(2))$ restrict to linear maps $\pa_e \colon \C A_1 \to \C A_{-1}$ and $\pa_f \colon \C A_{-1} \to \C A_1$. These restrictions can therefore be considered as densely defined unbounded operators at the level of Hilbert spaces, and upon doing so we denote them $\C E$ and $\C F$, respectively. Thus, $\C E \colon \Lambda(\C A_1)\to H_{-}$ and $\C F \colon \Lambda(\C A_{-1})\to H_{+}$, and one may show that $\C F \subseteq \C E^*$ and $\C E \subseteq \C F^*$. In particular, both $\C E$  and $\C F$ are closable and we denote their closures by $E$ and $F$ respectively. The D{\c a}browski-Sitarz Dirac operator is the unbounded selfadjoint operator on $H_{+}\oplus H_{-}$ with domain $\T{Dom}(E)\oplus \T{Dom}(F)$ given by
\[
D_q=\begin{pmatrix} 0 & F \\ E & 0 \end{pmatrix} .
\]
We denote the associated Lip-algebra  by
\[
C^{\Lip}(S_q^2):= \{ x\in C(S_q^2) \mid x(\T{Dom} (D_q))\subseteq \T{Dom}(D_q) \T{ and } [D_q,x] \T{ extends boundedly}\}.
\]
The map $\pa \colon  C^{\Lip}(S_q^2) \to \mathbb{B}(H_+\oplus H_- )$ given by $\pa(x)= \ov{[D_q,x]}$ is a derivation with respect to the diagonal action of $C^{\Lip}(S_q^2)$, and it has  the form
\[
\pa(x)= \begin{pmatrix} 0 & \pa_2(x) \\ \pa_1(x) & 0 \end{pmatrix} .
\]
The associated maps $\pa_1\colon C^{\Lip}(S_q^2) \to \mathbb{B}(H_+, H_-)$ and $\pa_2\colon C^{\Lip}(S_q^2) \to \mathbb{B}(H_-,H_+)$ are therefore also derivations for the natural $C^{\Lip}(S_q^2)$-bimodule structure on the two spaces, and on the coordinate algebra $\C O(S_q^2)$ one has the relations
\[
\pa_1(x)=q^{1/2} \rho(\pa_e(x))\vert_{H_+} \quad \text{ and } \quad \pa_2(x)=q^{-1/2} \rho(\pa_f(x))\vert_{H_-} .
\]
We obtain two natural seminorms
\begin{align*}
L_{D_q} & \colon \C O(S_q^2) \longrightarrow [0,\infty) \\
L_{D_q}^{\max} & \colon C^{\Lip}(S_q^2)  \longrightarrow [0,\infty),
\end{align*}
both given by taking the operator norm of the bounded extension of the commutator $[D_q,x]$. The main result in \cite{AgKa:PSM} is that $C(S_q^2)$ becomes a compact quantum metric space for the seminorm $L_{D_q}^{\max}$, and since $\C O(S_q^2) \su C^{\Lip}(S_q^2)$ this implies that the same is the case for $L_{D_q}$, see e.g.~\cite[Theorem 1.8]{Rie:MSA}.

\subsection{The quantum Berezin transform and quantum fuzzy spheres}
The key to the convergence results in \cite{AKK:PodCon} is a quantum analogue of the classical Berezin transform, and the main results in the present paper will turn out to be a consequence of a careful analysis of its analytic properties. Before embarking on this analysis, we first briefly introduce the quantum Berezin transform at the algebraic level, following \cite[Section 3]{AKK:PodCon}. As with the above, this will be kept rather short, and the reader is referred to \cite{AKK:PodCon} for more details.\\
For each $N\in \nn$, one obtains a state $h_N\colon C(S_q^2)\to \cc$ by setting
\begin{align}\label{eq:def-af-hN}
h_N(x) := \inn{N+1} \cd h\big( (a^*)^N x a^N \big), 
\end{align}
where $\inn{N+1}$ denotes the quantity $\sum_{k=0}^N q^{2k}$. The \emph{quantum Berezin transform} is then defined as the map $\be_N\colon C(S_q^2) \to C(S_q^2)$ given by
\[
\be_N(x):=(1\otimes h_N)\Delta(x).
\]
A priori, $\be_N(x)\in C(SU_q(2))$ but one can show that the image of $\be_N$ is actually contained in $C(S_q^2)$. Actually, it turns out that $\T{Im}(\be_N)$ is $(N+1)^2$-dimensional and contained in $\C O(S_q^2)$, and the quantum fuzzy sphere $F_q^N$ is defined as the concrete operator system $F_q^N:=\T{Im}(\be_N)$. In fact, it is proved in \cite[Lemma 3.3]{AKK:PodCon} that 
\[
F_q^N = \T{span}_{\cc}\{ u^{2n}_{in} \mid n \in \{0,1,\ldots,N\} \, , \, \, i \in \{0,1,\ldots,2n\} \} . 
\]

In the classical setting (i.e.~when $q=1$), one usually defines the fuzzy sphere in degree $N$ as the matrix algebra $\mathbb{M}_{N+1}(\cc)$, see \cite{Mad:Fuzz}. This is linked to the  2-sphere $S^2$ by means of the so-called covariant Berezin symbol $\sigma_N\colon \mathbb{M}_{N+1}(\cc) \to C(S^2)$ and its adjoint  
$\breve{\sigma}_N\colon C(S^2) \to \mathbb{M}_{N+1}(\cc)$  (see \cite[Section 2]{Rie:MSG} and references therein). 
The classical Berezin transform is then defined as $ \sigma_N \circ \breve{\sigma}_N $, and in \cite{AKK:PodCon} it was proven that this  agrees with the map $\beta_N$ constructed above. 
In the $q$-deformed setting, we have thus by-passed the covariant Berezin symbol and its adjoint and merely defined their composition, and the relationship between the quantum fuzzy sphere constructed above and its classical counterpart is given by the relation $F^N_1=\sigma_N\big(\mathbb{M}_{N+1}(\cc)\big)$.

\section{The Berezin transform is a Lip-norm contraction}
In \cite{AKK:PodCon} we showed that the quantum Berezin transform is a Lip-norm contraction on $\C{O}(S_q^2)$, and the aim of the current section is to prove that the same holds true at the level of 
$C^{\Lip}(S_q^2)$.
We consider again the orthonormal basis $\{ \Lambda(\ze^{n}_{ij}) \mid n \in \nn_0 \, , \, \, i,j \in \{0,1,\ldots,n\}\}$ for $L^2(SU_q(2))$ 
consisting of normalised matrix units, and define, for each $N \in \nn_0$, the linear map
$
\Phi_N \colon \B B(L^2(S_q^2)) \to F_q^N
$
by the formula
\[
\Phi_N(T) := \sum_{n = 0}^N \sum_{i = 0}^{2n} \ze^{2n}_{in} \inn{\Lambda(\ze^{2n}_{in}) ,T\La(1)} ,
\]
By construction, $\Phi_N$ is continuous from the weak operator topology (WOT) on  $\B B(L^2(S_q^2))$ to the norm topology on $F_{q}^N$.


\begin{lemma}\label{l:berphi}
For each $N \in \nn_0$ and $x \in C(S_q^2)$, it holds that $\be_N(x) = \be_N( \Phi_N(x) )$.
\end{lemma}
\begin{proof}
By norm-density of $\C O(S_q^2) \su C(S_q^2)$ and linearity and continuity of the involved operations, we only need to show that $\be_N(x) = \be_N(\Phi_N(x))$ in the case where $x = u^{2m}_{jm}$ for some $m \in \nn_0$ and $j \in \{0,1,\ldots,2m\}$. Suppose first that $m > N$. In this case we have that $\Phi_N(u^{2m}_{jm}) = 0$ since $\inn{\Lambda(u^{2k}_{ik}),\Lambda( u^{2m}_{jm})} = 0$ for all $k \leq N$ and all $i \in \{0,1,\ldots,2k\}$. Moreover, we see from  \cite[Lemma 3.2 and 3.3]{AKK:PodCon} 
that $\be_N(u^{2m}_{jm}) = 0$, hence proving the claimed identity in this case. Suppose next that $m \leq N$. We then have that
\begin{align}\label{eq:Phi-identity}
\Phi_N(u^{2m}_{jm}) = \sum_{k = 0}^N \sum_{i = 0}^{2k} \ze^{2k}_{ik} \inn{\Lambda(\ze^{2k}_{ik}) ,\Lambda(u^{2m}_{jm})}
= \ze^{2m}_{jm} \inn{\Lambda(\ze^{2m}_{jm}), \Lambda(u^{2m}_{jm})} = u^{2m}_{jm} ,
\end{align}
which proves the claimed identity in this case as well.
\end{proof}

We apply the notation $L^\infty(S_q^2)$ for the von Neumann algebra generated by $C(S_q^2) \su \B B(L^2(S_q^2))$, and define the extended Berezin transform $\wit{\be_N} \colon L^\infty(S_q^2) \to F_q^N$ as
\[
\wit{\be_N}(x) := \be_N( \Phi_N(x)) .
\]
Since $\Phi_N$ is WOT-norm continuous, we obtain that the same holds true for the extended Berezin transform $\wit{\be_N} \colon L^\infty(S_q^2) \to F_q^N$. Moreover, we see from Lemma \ref{l:berphi} that the extended Berezin transform agrees with the usual Berezin transform on $C(S_q^2) \su L^\infty(S_q^2)$.

\begin{lemma}\label{l:bercontra}
The extended Berezin transform $\wit{\be_N} \colon L^\infty(S_q^2) \to F_q^N$ is unital and completely positive and hence completely contractive.
\end{lemma}
\begin{proof}
The map $\beta_N\colon C(S_q^2) \to C(S_q^2)$ is defined by slicing the unital $\ast$-homomorphism $\De$ with the state $h_N$ and is therefore unital and completely positive. The fact that the same holds true for $\wit{\be_N}$ now follows by 
an approximation argument: given a positive element $x^*x\in L^\infty(S_q^2)$ there exists a net $\{x_i\}_{i \in I}$ in $C(S_q^2)$ converging to $x$ in the strong operator topology. The net $\{x_i^*x_i\}_{i \in I}$ then converges to $x^*x$ in the WOT and by the WOT-norm continuity of $\Phi_N$ we obtain that the net $\{ \beta_N(x_i^*x_i) \}_{i \in I}$ converges in norm to $\wit{\be_N}(x^*x)$. Since $\beta_N(x_i^*x_i)\geq 0$ for all $i \in I$, the positivity of $\wit{\be_N}(x^* x)$ follows. This shows that $\wit{\be_N}$ is positive and the fact that it is completely positive is proven in exactly the same manner.
\end{proof}

For each $j \in \{1,2,3\}$ we define the unbounded operator $\s D_j \colon \La\big( \C O(S_q^2)\big) \to L^2(S_q^2)$ by the formula $\s D_j( \La(x)) := \La( \de_j(x))$, where $\de_1, \de_2$ and $\de_3$ are the linear endomorphisms of $\C O(S_q^2)$ defined in \eqref{eq:delta-definitions}. Thus, $\s D_j$ is simply $\de_j$, but now thought of as an unbounded operator on $L^2(S_q^2)$. We also define $\s D_4 := -\s D_3 \colon \La\big( \C O(S_q^2)\big) \to L^2(S_q^2)$.

\begin{lemma}\label{l:adjointsdelta}
Each $\s D_j$, $j\in \{1,2,3,4\}$, is adjointable with  $\Lambda (\C O(S_q^2))\subseteq \T{Dom}(\s D_j^*)$ and for $y \in \C O(S_q^2)$ we have the explicit formulae:
\[
\begin{split}
\s D_1^*( \La(y)) & = q^{-1} \cd \s D_2( \La(y)),  \quad \s D_2^*( \La(y)) = q \cd \s D_1( \La(y)) \quad \mbox{and} \\
\s D_3^*( \La(y)) & = \s D_3( \La(y)), \quad \s D_4^*(\La(y)) = \s D_4(\La(y)).
\end{split}
\]
\end{lemma}
\begin{proof}
We only give the proof in the case of $\s D_1 \colon \La\big( \C O(S_q^2)\big) \to L^2(S_q^2)$ since the remaining cases follow from a similar argumentation. For $x,y \in \C O(S_q^2)$, we apply \eqref{eq:twisted-leibniz}, \eqref{eq:starder} and \eqref{eq:h-circ-delta=0} together with the defining relations in the Hopf $*$-algebra $\C U_q( \G{su}(2))$ (see \cite[Page 4]{AKK:PodCon}) to compute as follows:
\begin{align*}
\inn{ \La(y), \s D_1 \La(x)}  & = h\big(y^* \de_1(x)\big) = h\big( \de_1( \de_k(y^*) x)\big) - h\big( \de_1(\de_k(y^*)) \de_k(x)\big)  \\
& = -  h\big( \de_1(\de_k(y^*)) \de_k(x)\big) = - q^{-1} h\big( \de_k( \de_1(y^*)) \de_k(x)\big) \\
& = - q^{-1} h( \de_1(y^*) x) = q^{-1} h( \de_2(y)^* x) = \inn{q^{-1} \cd \s D_2 \La(y), \La(x)} . \qedhere
\end{align*}
\end{proof}

For each $N \in \nn_0$ we consider the orthogonal projection $P_N \colon L^2(S_q^2) \to L^2(S_q^2)$ given by the formula
\[
P_N(\xi) := \sum_{n = 0}^N \sum_{i = 0}^{2n} \Lambda(\ze^{2n}_{in}) \inn{\Lambda(\ze^{2n}_{in}) ,\xi } .
\]
Remark that $P_N \La(x) = \La( \Phi_N(x))$ for all $x \in C(S_q^2)$.

\begin{lemma}\label{l:polycommut}
For each $j \in \{1,2,3,4\}$ and each $N \in \nn_0$, it holds that
\[
P_N \s D_j \su \s D_j P_N.
\]
\end{lemma}
\begin{proof}
Let $N \in \nn_0$ be given. We only need to consider the case where $j \in \{1,2,3\}$. Using the description of $\de_j : \C O(S_q^2) \to \C O(S_q^2)$ in terms of the coproduct on $\C O(SU_q(2))$ and the pairing of Hopf $*$-algebras $\inn{\cd,\cd} : \C U_q(\G{su}(2)) \ti \C O(SU_q(2)) \to \cc$ we obtain that 
\[
\begin{split}
\Phi_N \de_j( u^{2n}_{in} ) 
& = \fork{ccc}{ 0 & \T{for} & n > N \\ \de_j(u^{2n}_{in}) & \T{for} & 0 \leq n \leq N } \\
& = \de_j \Phi_N(u^{2n}_{in}) .
\end{split}
\]
The above identity now proves the lemma since $\T{Dom}(\s D_j)$ agrees with the linear span of the elements $\La(u^{2n}_{in} ) \in L^2(S_q^2)$ (for $n \in \nn_0$ and $i \in \{0,1,\ldots,2n\}$).
%
\end{proof}

We introduce the linear map
\[
\de \colon \C O(S_q^2) \to M_2\big( \C O(S_q^2) \big) \quad \de(x) := \ma{cc}{-\de_3(x) & \de_2(x) \\ \de_1(x) & \de_3(x)} .
\]
An application of \cite[Proposition 3.11]{AKK:PodCon} then shows that
\[
\de(x) = u \pa(x) u^* \quad \T{for all } x \in \C O(S_q^2),
\]
where $u$ is the fundamental unitary corepresentation. We then introduce the four linear maps $\wit{\de_j} \colon C^{\Lip}(S_q^2) \to \B B(L^2(S_q^2))$, $j \in \{1,2,3,4\}$, by the formulae
\[
\begin{split}
\wit{\de_1}(x) & := (u \pa(x) u^*)_{10},  \quad \wit{\de_2}(x) := (u \pa(x) u^*)_{01} \quad \T{and} \\
\wit{\de_3}(x) & := (u \pa(x) u^*)_{11}, \quad \wit{\de_4}(x) := (u \pa(x) u^*)_{00}, 
\end{split}
\]
where the subscripts denote the matrix entries of $u \pa(x) u^*\in \mathbb{M}_2\big(\mathbb{B}(L^2(S_q^2))\big)$. We also introduce the linear map
\[
\wit{\de} \colon C^{\Lip}(S_q^2) \to M_2\big( \B B(L^2(S_q^2)) \big) \quad 
\wit{\de}(x) := \ma{cc}{ \wit{\de_4}(x) & \wit{\de_2}(x) \\ \wit{\de_1}(x) & \wit{\de_3}(x)} .
\]
By construction we have that $\wit{\de}(x) = u \pa(x) u^*$ for all $x \in C^{\Lip}(S_q^2)$ and hence from the above we obtain that $\wit{\de}(x) = \de(x)$ for all $x \in \C O(S_q^2)$.

We wish to consider each $\wit{\de_j}$ as a densely defined unbounded operator. More precisely, for 
each $j \in \{1,2,3,4\}$ we define an unbounded operator $\wit{\s D_j} \colon \La\big(C^{\Lip}(S_q^2) \big) \to L^2(S_q^2)$ by setting
\[
\wit{\s D_j}( \La(x)) := \wit{\de_j}(x)\big(  \La(1)  \big).
\]
By construction, the unbounded operator $\wit{\s D_j}$ is an extension of the unbounded operator $\s D_j \colon \La( \C O(S_q^2)) \to L^2(S_q^2)$. 
We are now going to show that $\wit{\s D_j}$ also has a densely defined adjoint. To do so, it turns out to be convenient to work with a slightly more general type of unbounded operator. For each $\xi,\eta \in \C A_1$ (see \eqref{eq:def-of-spectral-subspaces} for a definition of this space) we define the unbounded operators $R_{\xi,\eta}, T_{\eta^*,\xi^*}\colon \La\big(C^{\Lip}(S_q^2)\big) \to L^2(S_q^2)$ by
\[
R_{\xi,\eta}( \La(x)) := (\xi \cd \pa_1(x) \cd \eta)  \La(1) \quad \T{and} \quad T_{\eta^*,\xi^*}(\La(x)) := (\eta^* \cd \pa_2(x) \cd \xi^*) \La(1) ,
\]
where we consider $\xi$ as a bounded multiplication operator from $H_-$ to $H_0$, $\eta$ as a bounded multiplication operator from $H_0$ to $H_+$ and $\pa_1(x)$ as a bounded operator from $H_+$ to $H_-$ and $\pa_2(x)$ as a bounded operator from $H_-$ to $H_+$. 

\begin{lemma}
For $\xi,\eta \in \C A_1$,  the unbounded operators $R_{\xi,\eta}$ and $T_{\eta^*,\xi^*}$ are adjointable with $\La\big( \C O(S_q^2) \big)$ in the domain of their adjoints, and on $\La\big( \C O(S_q^2) \big)$ we have the explicit formulae
\[
\begin{split}
R_{\xi,\eta}^*\La(y) 
& = \La\big( \pa_f( \xi^*y) \cd \nu(\eta)^*
- \xi^* y \cd \nu(\pa_e(\eta))^* \big) \quad \mbox{and} \\
T_{\eta^*,\xi^*}^* \La(y) & = \La\big( \pa_e( \eta y) \cd \nu(\xi^*)^* - \eta y \cd \nu(\pa_f(\xi^*))^* \big) 
\end{split}
\] 
for all $y \in \C O(S_q^2)$. In particular, we obtain that $R_{\xi,\eta}^* \La(y)$ and $T_{\eta^*,\xi^*}^* \La(y)$ belong to $\La\big( \C O(S_q^2)\big)$.
\end{lemma}
\begin{proof}
We only give a proof in the case of $R_{\xi,\eta}$ since the proof for $T_{\eta^*,\xi^*}$ follows a similar structure. Let $x \in C^{\Lip}(S_q^2)$ and $y \in \C O(S_q^2)$ be given. Using that $\pa_1(x)=[E, x]$ on $\T{Dom}(E)\supseteq \Lambda(\C A_1)$ in combination with the twisted tracial property of the Haar state from \eqref{eq:twisttrace} we compute as follows:
\[
\begin{split}
\inn{\La(y), R_{\xi,\eta}( \La(x)) } & = \binn{\La(y), \xi \cd \pa_1(x)( \La(\eta) )} \\
& = \binn{\La(\xi^* y), E(x \La(\eta) )} - \binn{\La(\xi^*y), x  E( \La(\eta))} \\
& = \binn{F\big( \La(\xi^* y) \big), \La( x \cd \eta ) } - \binn{\La(\xi^*y), \La( x \pa_e(\eta))} \\
& = \binn{ \La\big( \pa_f(\xi^* y) \cd \nu(\eta)^*\big), \La(x)} -\binn{\La\big(\xi^*y \cd \nu(\pa_e(\eta))^* \big), \La(x)} .
\end{split}
\]
This proves the result of the lemma.
\end{proof}

Since $\wit{\de}(x)= u \pa(x) u^*$, a direct computation of the matrix products shows that each entry is of the form $R_{\xi_1,\eta_1}+T_{\eta_2^*,\xi_2^*}$ for some $\xi_1,\eta_1, \xi_2, \eta_2\in \C A_1$, and hence it follows from the lemma above that each of the corresponding  unbounded operators  $\wit{ \s D_j} : \La( C^{\Lip}(S_q^2)) \to L^2(S_q^2)$ has a densely defined adjoint and that
\[
\wit{\s D_j}^*\big( \La( \C O(S_q^2)) \big) \su \La( \C O(S_q^2)) .
\]
We now prove that the conclusion of Lemma \ref{l:polycommut} also holds true for the extended unbounded operators.
\begin{lemma}\label{l:phidel}
For each $N \in \nn_0$ and $j \in \{1,2,3,4\}$ we have that
\[
P_N \wit{\s D_j} \su \s D_j P_N .
\]
\end{lemma}
\begin{proof}
Let $x \in C^{\Lip}(S_q^2)$ be given. It suffices to show that
\[
\inn{P_N \wit{\s D_j}\big( \La(x) \big), \xi } = \inn{\s D_j P_N(\La(x)), \xi}
\]
for all $\xi \in L^2(S_q^2)$. 
 Since $\wit{\s D_j}^*( \La(\C O(S_q^2))) \su \La(\C O(S_q^2))$ there exists an $M \geq N$ such that $P_M\wit{\s D_j}^* P_N \xi = \wit{\s D_j}^* P_N \xi$. Using Lemma \ref{l:polycommut} and the inclusion $\s D_j \su \wit{\s D_j}$ we may thus compute as follows:
\[
\begin{split}
\inn{P_N \wit{\s D_j}\big( \La(x) \big), \xi} & = \inn{ \La(x), P_M \wit{\s D_j}^* P_N \xi} 
= \inn{ P_N \wit{\s D_j} P_M \La(x), \xi} \\
& = \inn{ P_N \s D_j P_M \La(x), \xi} 
= \inn{ \s D_j P_N P_M \La(x), \xi} = \inn{\s D_j P_N \La(x), \xi} . \qedhere
\end{split}
\]
\end{proof}

To ease the notation, we define an algebra automorphism $\nu^{1/2} \colon \C O(SU_q(2)) \to \C O(SU_q(2))$ as the composition $\nu^{1/2} := \pa_k^{-1} \ci \de_k^{-1}$. Notice that the square root notation indeed makes sense since $\nu^{1/2} \ci \nu^{1/2} = \nu$. We furthermore define an antilinear surjective isometry $J \colon L^2(SU_q(2)) \to L^2(SU_q(2))$ by the formula 
\[
J\big( \La(y) \big) := \La\big( \nu^{-1/2}(y^*)\big) \q \T{for all } y \in \C O(SU_q(2)) .
\]
Clearly $J^2 = \T{id}$ and using that $\nu^{\frac12}(x)^*=\nu^{-\frac12}(x^*)$  (see \cite[(2.10) \& (3.3)]{AKK:PodCon}), a direct computation shows that
\begin{equation}\label{eq:tomtak}
[ J \rho(x) J, \rho(y)] = 0 \quad \T{ for all } x,y \in C(SU_q(2)),
\end{equation}
where $\rho \colon C(SU_q(2)) \to \B B(L^2(SU_q(2)))$ as usual denotes the GNS-representation. In fact, it can be verified that $J \colon L^2(SU_q(2)) \to L^2(SU_q(2))$ is the phase operator in the polar decomposition of the adjoint operation $* \colon \La( \C O(SU_q(2))) \to L^2(SU_q(2))$ (considered as an anti-linear unbounded operator on the GNS-space). The analysis we are carrying out is therefore fully compatible with Tomita-Takesaki theory and $J$ is exactly the modular conjugation associated to $\C O(SU_q(2))$ considered as a left Hilbert algebra with inner product coming from the Haar-state, see \cite[Chapter VI]{Tak:TOA}. 
To continue, we notice that $J( H_{\pm}) = H_{\mp}$ and we apply the same notation $J$ for the induced map
\[
J := \ma{cc}{0 & J \\ -J & 0} : H_+ \op H_- \longrightarrow H_+ \op H_- .
\] 
The latter antilinear surjective isometry $J \colon H_+ \op H_- \to H_+ \op H_-$ is referred to as the reality operator since it provides the even unital spectral triple $(C(S_q^2),H_+ \op H_-, D_q)$ with a real structure of dimension $2$, see \cite[Definition 3]{Con:NGR}. This real structure was also found by D{\c a}browski and Sitarz in their paper \cite{DaSi:DSP}. For the convenience of the reader we provide some details on these matters. We first of all have the identities $J^2 = -\T{id}$ and $J \ga = -\ga J$ (here $\gamma$ denotes the grading operator) at the level of operations on $H_+ \op H_-$ and it can be verified that 
\begin{equation}\label{eq:jaycommut}
J( \T{Dom}(D_q)) = \T{Dom}(D_q) \quad \T{and} \quad D_q J = J D_q .
\end{equation}
Secondly, to prove the first order condition, i.e.~that $[ JyJ^{-1},\pa(x) ] = 0$ for all $x,y \in C^{\Lip}(S_q^2)$,  we may, without loss of generality, assume that $y \in \C O(S_q^2)$. Letting $\xi \in \T{Dom}(D_q)$ and applying \eqref{eq:tomtak} and \eqref{eq:jaycommut} we then compute that
\[
\begin{split}
[J y J^{-1}, \pa(x)](\xi) 
& = [J y J^{-1}, D_q x - x D_q](\xi)
= [J y J^{-1}, D_q] x(\xi) - x [JyJ^{-1},D_q](\xi)  \\
& = \big[x, [D_q, JyJ^{-1}] \big](\xi) = [x, J \pa(y) J^{-1}](\xi) .
\end{split}
\]
Since $\pa_1(y), \pa_2(y) \in \C O(SU_q(2))$ and $J \pa(y) J^{-1} = \ma{cc}{ 0 &  -J\pa_1(y)J \\ -J \pa_2(y) J & 0}$ we now immediately obtain from \eqref{eq:tomtak} that $[J y J^{-1},\pa(x)] = 0$. \\

We are now ready to prove that the maps $\wit{\de_j}$ take values in $L^\infty(S_q^2)$. This is essential since it allows us to compose these maps with the extended Berezin transform. In fact, we shall see that $\wit{\de_j}$ actually commutes with the extended Berezin transform and this will be the key ingredient in our proof of Theorem \ref{mainthm:B}.

\begin{lemma}\label{lem:delta-in-L-infty}
For $j \in \{1,2,3,4\}$ and $x \in C^{\Lip}(S_q^2)$, it holds that $\wit{\de_j}(x) \in L^\infty(S_q^2)$.
\end{lemma}
\begin{proof}
We first note that the antilinear surjective isometry $J \colon L^2(SU_q(2)) \to L^2(SU_q(2))$ restricts to an antilinear surjective isometry $J \colon L^2(S_q^2) \to L^2(S_q^2)$. As above this restricted $J$ agrees with the modular conjugation arising from $\C O(S_q^2)$ considered as a left Hilbert algebra with inner product coming from the Haar state $h \colon C(S_q^2) \to \cc$. In particular, it follows from \cite[Theorem 1.19]{Tak:TOA} that the commutant of $L^\infty(S_q^2) \su \B B(L^2(S_q^2))$ agrees with $J L^\infty(S_q^2) J$. By von Neumann's bicommutant theorem, it hence suffices to show that 
\[
[ J y J^{-1},  u \pa(x) u^*] = 0
\]
for all $y \in \C O(S_q^2)$. This does however follow from the first order condition $[J y J^{-1}, \pa(x)] = 0$ and the fact that the fundamental corepresentation unitary $u$ belongs to $\mathbb{M}_2\big( \C O(SU_q(2))\big)$. Indeed, using \eqref{eq:tomtak} we obtain that 
\[
[ J y J^{-1},  u \pa(x) u^*] =u[JyJ^{-1}, \pa(x)]u^* = 0 . \qedhere
\]
\end{proof}

We also record the following easy consequence of the above lemmas:

\begin{cor}\label{c:phidel}
For each $N \in \nn_0$ and each $j \in \{1,2,3,4\}$  we have $\de_j \Phi_N(x)=\Phi_N\wit{\de_j}(x)$ for all  $x\in C^{\Lip}(S_q^2)$.
\end{cor}
\begin{proof}
The map $\Lambda \colon C(S_q^2)\to L^2(S_q^2)$ extends (injectively) to $L^\infty(S_q^2)$ and it is straightforward to see that the relation $P_N\Lambda(z)= \Lambda \Phi_N(z)$ still holds for $z\in L^\infty(S_q^2)$. For $x\in C^{\Lip}(S_q^2)$ we therefore obtain (using  Lemma \ref{l:phidel}) that
\begin{align*}
\Lambda(\de_j\Phi_N(x))&=\s D_j(\Lambda(\Phi_N(x)))
=\s D_j P_N\Lambda(x)
= P_N \wit{\s D_j} \Lambda (x) \\
& =P_N(\Lambda(\wit{\de_j} (x)))
=\Lambda(\Phi_N\wit{\de_j} (x)) . \qedhere
\end{align*}
\end{proof}

\begin{prop}\label{p:delber}
For each $N \in \nn_0$ and $j \in \{1,2,3,4\}$ we have that $\de_j \be_N(x) = \wit{\be_N} \wit{\de_j}(x)$ for all $x \in C^{\Lip}(S_q^2)$.
\end{prop}
\begin{proof}
First note that the composition on the right hands side indeed makes sense by Lemma \ref{lem:delta-in-L-infty}.
Let $x \in C^{\Lip}(S_q^2)$ be given. By \cite[Lemma 3.7]{AKK:PodCon}, 
 Lemma \ref{l:berphi} and Corollary \ref{c:phidel} we obtain that
\[
\de_j \be_N(x) = \de_j \be_N \Phi_N(x) = \be_N \de_j \Phi_N(x) = \be_N \Phi_N \wit{\de_j}(x) = \wit{\be_N} \wit{\de_j}(x) . \qedhere
\]
\end{proof}

With the above results at our disposal, it is now an easy task to show that the Berezin transform is a Lip-norm contraction also at the level of the Lipschitz algebra $C^{\Lip}(S_q^2)$. This extends the result in \cite[Proposition 3.12]{AKK:PodCon}.

\begin{theorem}\label{t:berlipcon}
For each $N \in \nn_0$ and $x \in C^{\Lip}(S_q^2)$ we have the inequality
\[
L_{D_q}^{\max}(\be_N(x)) \leq L_{D_q}^{\max}(x) .
\]
\end{theorem}
\begin{proof}
By \cite[Proposition 3.11]{AKK:PodCon}, Lemma \ref{l:bercontra} and Proposition \ref{p:delber} it holds that
\[
L_{D_q}^{\max}(\be_N(x)) = \| u \pa( \be_N(x)) u^* \| = \| \de( \be_N(x)) \| = \| \wit{\be_N}( \wit{\de}(x)) \| \leq \| \wit{\de}(x) \| = L_{D_q}^{\max}(x). \qedhere
\]
\end{proof}

\section{The coproduct commutes with the Dirac operator}
In order to prove our main results, we are still missing one important ingredient, namely a certain compatibility between the Dirac operator and the coproduct. This compatibility result is related to the $\C U_q(\G{su}(2))$-equivariance of the Dirac operator which was described in \cite{DaSi:DSP}. We are here clarifying that this equivariance actually comes from an underlying relationship between the coproduct and the Dirac operator. In order to achieve this, we need the unitary operator $W \colon L^2(SU_q(2)) \hot L^2(SU_q(2)) \to L^2(SU_q(2)) \hot L^2(SU_q(2))$ defined by the formula
\[
W( \La(x) \ot \La(y)) := \De(y)  \big( \La(x) \ot \La(1) \big) .
\]
This operator (known as the \emph{multiplicative unitary} for $SU_q(2)$) implements the coproduct $\De \colon C(SU_q(2)) \to C(SU_q(2)) \ot_{\T{min}} C(SU_q(2))$ in the sense that
\[
\De(z) = W(1 \ot z) W^* \quad \T{for all } z \in C(SU_q(2)) ,
\]
see \cite{BaSk:UMD} for more details on multiplicative unitaries. This formula also shows that $\De$ extends to a normal $*$-homomorphism 
\[
\De\colon L^\infty(SU_q(2)) \to  L^\infty(SU_q(2)) \bar{\otimes}  L^\infty(SU_q(2)),
\]
and this actually turns $ L^\infty(SU_q(2))$ into a von Neumann algebraic compact quantum group. 
Using \cite[Theorem 1.2]{Wor:UAC}, it is not difficult to see that $\C A_1$ is generated by $a$ and $b$ as a right module over $\C O(S_q^2)$, and, similarly, that $\C A_{-1}$ is generated by $a^*$ and $b^*$.  Since
\[
\Delta(a)=a\ot a - qb^*\ot b \quad \T{ and } \quad \Delta(b)=b\ot a + a^*\ot b,
\]
it therefore follows that   $\C A_1$ and $\C A_{-1}$ are comodules for $\C O(SU_q(2))$ in the sense that $\Delta (\C A_i)\subseteq \C O(SU_q(2))\otimes \C A_i$, $i=\pm 1$, and the
 unitary operator $W$  therefore restricts to two unitary  operators
\[
\begin{split}
W_+ & \colon L^2(SU_q(2)) \hot H_+  \longrightarrow L^2(SU_q(2)) \hot H_+ \quad \T{and} \\
W_- & \colon L^2(SU_q(2)) \hot H_-  \longrightarrow L^2(SU_q(2)) \hot H_- .
\end{split}
\]
We now define the unbounded selfadjoint operator 
\[
1 \hot D_q \colon \T{Dom}(1 \hot D_q) \to L^2(SU_q(2)) \hot (H_+ \op H_-)
\]
as the closure of the unbounded symmetric operator
\[
1 \ot \ma{cc}{0 & \C F \\ \C E & 0} \colon \La\big( \C O(SU_q(2)) \big) \ot \big( \La(\C A_1) \op \La(\C A_{-1}) \big) \longrightarrow L^2(SU_q(2)) \hot (H_+ \op H_-) .
\] 

\begin{lemma}\label{l:multiunidir}
The unitary operator $W_+ \op W_-$ preserves $\T{Dom}(1 \hot D_q)$ and
\[
[ 1 \hot D_q , W_+ \op W_- ](\xi) = 0
\]
for all $\xi \in \T{Dom}(1 \hot D_q)$.
\end{lemma}
\begin{proof}
Since $W_+ \op W_-$ preserves the core $\La\big( \C O(SU_q(2)) \big) \ot \big( \La(\C A_1) \op \La(\C A_{-1}) \big)$, it suffices to  verify the identity for 
\[
\xi = \La(x) \ot \ma{c}{\La(y_1) \\ \La(y_{-1})}
\]
with $x \in \C O(SU_q(2))$ and $y_1 \in \C A_1$ and $y_{-1} \in \C A_{-1}$. It then follows that $W_+ \op W_-$ preserves $\T{Dom}(1 \hot D_q)$ and that the commutation relation holds true here as well. In this case, the desired identity follows from \cite[Lemma 4.1]{AKK:PodCon} via the computation:
\begin{align*}
 (1 \hot D_q) (W_+ \op W_-)(\xi) 
&  = \ma{c}{(1 \ot \C F)\big( \De(y_{-1})  ( \La(x) \ot \La(1) ) \big) \\ (1 \ot \C E)\big( \De(y_{1}) ( \La(x) \ot \La(1) ) \big)} \\
&  = \ma{c}{\De( \pa_f(y_{-1}) )  ( \La(x) \ot \La(1) )  \\  \De(\pa_e(y_{1}))  ( \La(x) \ot \La(1) ) } \\
&= (W_+ \op W_-)(1 \hot D_q)(\xi) . \qedhere
\end{align*}
\end{proof}

The unital $C^*$-algebra $C(SU_q(2)) \ot_{\T{min}} C(S_q^2)$ acts on $L^2(SU_q(2)) \hot (H_+ \op H_-)$ via the representation $\rho \ot \pi$ (which we will from now on often suppress). We define the dense unital $*$-subalgebra  
\[
\Lip_{1 \hot D_q}(C(SU_q(2)) \ot_{\T{min}} C(S_q^2)) \su C(SU_q(2)) \ot_{\T{min}} C(S_q^2)
\]
to consist of those $x \in  C(SU_q(2)) \ot_{\T{min}} C(S_q^2)$ which preserves the domain of $1\hot D_q$ and whose commutator with $1\hot D_q$ extends to a bounded operator on $L^2(SU_q(2)) \hot ( H_+ \op H_-) $.
For each $x \in \Lip_{1 \hot D_q}\left(C(SU_q(2)) \ot_{\T{min}} C(S_q^2)) \right)$ we let
\[
(1 \ot \pa)(x) \colon L^2(SU_q(2)) \hot ( H_+ \op H_-) \longrightarrow L^2(SU_q(2)) \hot ( H_+ \op H_-)
\]
denote the bounded extension of the commutator 
\[
[1 \hot D_q,x] : \T{Dom}(1 \hot D_q) \longrightarrow L^2(SU_q(2)) \hot ( H_+ \op H_-) .
\]
The next result shows that $\pa$ commutes with the comultiplication, thus providing us with an analytic generalisation of \cite[Lemma 4.1]{AKK:PodCon}.

For each $x \in C^{\Lip}(S_q^2)$ we define the bounded operator
\[
\begin{split}
\De(\pa(x)) & := (W_+ \op W_-)( 1 \ot \pa(x) )(W_+ \op W_-)^* \\
& \quad \quad \colon L^2(SU_q(2)) \hot (H_+ \op H_-) \to L^2(SU_q(2)) \hot (H_+ \op H_-) .
\end{split}
\]

\begin{lemma}\label{l:partialcommutes}
For $x \in C^{\Lip}(S_q^2)$ it holds that $\De(x) \in \Lip_{1 \hot D_q}\big(C(SU_q(2)) \ot_{\T{min}} C(S_q^2) \big)$ and
$(1 \ot \pa)\De(x) = \De(\pa(x))$.  
\end{lemma}
\begin{proof}
We first remark that $(i + (1\hot D_q) )^{-1}= 1 \ot (i+D_q)^{-1}$ so that
\[
\Lambda\big(\C O(SU_q(2))\big)\ot \T{Dom}(D_q)= (i + (1\hot D_q) )^{-1} \Big(\Lambda\big(\C O(SU_q(2))\big) \ot L^2(S_q^2) \Big).
\]
From this it follows that $\Lambda\big(\C O(SU_q(2))\big)\ot \T{Dom}(D_q)$ is a core for $1\hot D_q$. We then obtain that
$1 \ot x \in \Lip_{1 \hot D_q}\big(C(SU_q(2)) \ot_{\T{min}} C(S_q^2)\big)$ and moreover that $(1 \ot \pa)(1 \ot x) = 1 \ot \pa(x)$. Let now $\xi \in \T{Dom}(1 \hot D_q)$. Using Lemma \ref{l:multiunidir} we see that
\[
\De(x)(\xi) = (W_+ \op W_-) (1 \ot x) (W_+^* \op W_-^*)(\xi) \in \T{Dom}(1 \hot D_q)
\]
and compute the commutator
\[
\begin{split}
[ (1 \hot D_q), \De(x) ](\xi) 
& = (W_+ \op W_-) [ (1 \hot D_q) , (1 \ot x) ] (W_+^* \op W_-^*)(\xi) \\
& = (W_+ \op W_-)(1 \ot \pa(x) )  (W_+^* \op W_-^*)(\xi) = \De(\pa(x))(\xi) .
\end{split}
\]
This proves the lemma.
\end{proof}

For each $\xi,\ze \in L^2(SU_q(2))$ we let $\phi_{\xi,\ze} \colon C(SU_q(2)) \to \cc$ denote the bounded linear functional $\phi_{\xi,\ze}(x) := \inn{\xi, \rho(x) \ze}$.

\begin{lemma}\label{lem:slice-and-Lip}
For each $\xi,\ze \in L^2(SU_q(2))$ and $z \in \Lip_{1 \hot D_q}\big(C(SU_q(2)) \ot_{\T{min}} C(S_q^2)\big)$ it holds that
$(\phi_{\xi,\ze} \ot 1)(z) \in  C^{\Lip}(S_q^2)$ and that
\[
\pa\big( (\phi_{\xi,\ze} \ot 1)(z) \big) = (\phi_{\xi,\ze} \ot 1)(1 \ot \pa)(z) .
\]
In particular, we have the estimate
\[
L_{D_q}^{\max}\big( (\phi_{\xi,\ze} \ot 1)(\De(x)) \big) \leq \| \xi \| \| \ze\| \cd L_{D_q}^{\max}(x)
\]
for all $x \in C^{\Lip}(S_q^2)$.
\end{lemma}
\begin{proof}
We define two bounded operators $T_\xi, T_\ze \colon H_+ \op H_- \to L^2(SU_q(2)) \hot (H_+ \op H_-)$ by $T_\xi(\eta) := \xi \ot \eta$ and $T_\ze(\eta) := \ze \ot \eta$. It can then be proved that $T_\ze(\T{Dom}(D_q)) \su \T{Dom}(1 \hot D_q)$ and $T_\xi^*\big( \T{Dom}(1 \hot D_q)\big) \su \T{Dom}(D_q)$. Moreover, it holds that 
\[
T_\ze D_q \su (1 \hot D_q) T_\ze  \q \T{and} \q T_\xi^* (1 \hot D_q) \su D_q T_\xi^* .
\]
We thus obtain that $(\phi_{\xi,\ze} \ot 1)(z) = T_\xi^* z T_\ze$ preserves the domain of $D_q$ and that
\[
[ D_q, T_\xi^* z T_\ze ] \su T_\xi^* [1 \hot D_q, z] T_\ze \su T_\xi^* (1 \ot \pa)(z) T_\ze = (\phi_{\xi,\ze} \ot 1)(1 \ot \pa)(z) .
\]
This proves the first part of the lemma. The second part now follows immediately by an application of Lemma \ref{l:partialcommutes}.
\end{proof}

\section{Proofs of the main results}
In this section we gather the proofs of our main results, which are now easily obtained with the tools developed in the previous sections at our disposal. Recall that the co-amenability of $SU_q(2)$ means that the counit $\epsilon\colon \C O (SU_q(2))\to\cc$ extends to a $*$-character on $C(SU_q(2))$, and by restricting its domain we may consider $\epsilon$ as an element in $\C S (C(S_q^2)))$.  We also have the sequence of states $\{h_N\}_{N = 1}^\infty$ in $\C S (C(S_q^2))$ defined in \eqref{eq:def-af-hN} and hence we may consider the distance $d_q^{\max}(\epsilon, h_N)$ with respect to the metric on $\C S (C(S_q^2))$ defined via the maximal seminorm $L_{D_q}^{\max}\colon C^{\Lip}(S_q^2)\to [0,\infty)$. The following proposition now shows that the quantum Berezin transforms approximate the identity map on the Lip-unit ball well:

\begin{prop}\label{prop:berezin-approximates-identity}
For each $x\in C^{\Lip}(S_q^2)$ it holds that $\|x-\be_N(x)\|\leq d_q^{\max}(h_N,\epsilon)L_{D_q}^{\max}(x)$
\end{prop}

In \cite[Proposition 4.3]{AKK:PodCon} the corresponding statement was proven for $d_q$ (the metric arising from restricting the domain of $L_{D_q}$ to $\C O(S_q^2)$) and $x\in \C O (S_q^2)$, but with the above analysis at our disposal the proof carries over verbatim:

\begin{proof}
First note that $\be_N(x)-x=(1\otimes (h_N-\epsilon))\Delta(x)$, so for unit vectors $\xi,\eta\in L^2(SU_q(2))$ we have
\begin{align*}
\phi_{\xi, \eta}(\be_N(x)-x)= (h_N-\epsilon) ((\phi_{\xi,\eta} \otimes 1) \Delta(x))\leq d_q^{\max}(h_N, \epsilon) L_{D_q}^{\max}(x),
\end{align*}
where the last inequality follows from Lemma \ref{lem:slice-and-Lip}.
\end{proof}

The estimates in Proposition \ref{prop:berezin-approximates-identity} and Theorem \ref{t:berlipcon} together with the results from  \cite{AgKa:PSM} and \cite{AKK:PodCon}, now allow us to prove our main results:

\begin{proof}[Proof of Theorem \ref{mainthm:B}] 
By \cite[Proposition 4.4]{AKK:PodCon}, we know that $h_N$ converges to $\epsilon$ in the weak$^*$-topology on $\C S(C(S_q^2))$ and by \cite[Theorem 8.3]{AgKa:PSM} it holds that $C(S_q^2)$ is a compact quantum metric space with respect to $L_{D_q}^{\max}$. We may thus conclude that $\lim_{N\to \infty}d_q^{\max}(h_N, \epsilon)=0$. Theorem \ref{mainthm:B} now follows immediately from the estimates in Proposition \ref{prop:berezin-approximates-identity} and Theorem \ref{t:berlipcon} since $\be_N(x)\in \C O(S_q^2)$ for all $N\in \nn$ and all $x\in C^{\Lip}(S_q^2)$
\end{proof}

\begin{proof}[Proof of Theorem \ref{mainthm:A}]
This follows immediately from Theorem \ref{mainthm:B} and \cite[Lemma 4.15]{AKK:PodCon}.
\end{proof}


\end{document}